\documentclass{amsart}
\usepackage{amssymb, amsthm, amsmath, amscd}

\newtheorem{theorem}{Theorem}[section]
\newtheorem{lemma}[theorem]{Lemma}
\newtheorem{prop}[theorem]{Proposition}
\newtheorem{cor}[theorem]{Corollary}

\theoremstyle{definition}
\newtheorem{defn}[theorem]{Definition}
\newtheorem{example}[theorem]{Example}

\theoremstyle{remark}
\newtheorem{remark}[theorem]{Remark}

\def\nil{\mathrm{nil}}
\def\sep{\mathrm{sep}}
\def\st{\mathrm{st}}

\def\ann{\operatorname{ann}}
\def\rad{\operatorname{rad}}
\def\Ext{\operatorname{Ext}}
\def\Proj{\operatorname{Proj}\,}
\def\Spec{\operatorname{Spec}}
\def\depth{\operatorname{depth}}
\def\image{\operatorname{image}}
\def\to{\longrightarrow}
\def\into{\hookrightarrow}
\def\hat{\widehat}

\def\a{\mathfrak a}
\def\b{\mathfrak b}
\def\p{\mathfrak p}
\def\q{\mathfrak q}
\def\m{\mathfrak m}
\def\n{\mathfrak n}

\def\A{\mathfrak A}
\def\B{\mathfrak B}
\def\P{\mathfrak P}
\def\Q{\mathfrak Q}
\def\M{\mathfrak M}

\def\FF{\mathbb F}
\def\NN{\mathbb N}

\def\QQ{\mathbb Q}
\def\ZZ{\mathbb Z}

\def\calO{\mathcal O}

\begin{document}
\title{A connectedness result in positive characteristic}

\author{Anurag K. Singh}
\address{Department of Mathematics, University of Utah, 155 South 1400 East,
Salt Lake City, UT~84112, USA} \email{singh@math.utah.edu}
\author{Uli Walther}
\address{Department of Mathematics, Purdue University, 150 N. University Street,
West Lafayette, IN~47907, USA} \email{walther@math.purdue.edu}

\dedicatory{Dedicated to Professor Paul Roberts on the occasion of his sixtieth
birthday}

\thanks{The first author was supported by the NSF under Grant DMS~0300600, and
the second author under Grant DMS~0100509. Both authors are grateful to Gennady
Lyubeznik for useful discussions and comments.}

\date{\today}
\maketitle

\section{Introduction}

All rings considered in this note are commutative and Noetherian. We give a
simple proof of the following result due to Lyubeznik:

\begin{theorem}\cite[Corollary~4.6]{Ly-dual}\label{connected}
Let $(R,\m)$ be a complete local ring of positive dimension, with a separably
closed coefficient field of positive characteristic. Then the $e$-th iteration
of the Frobenius map
$$
F:H^1_\m(R)\to H^1_\m(R)
$$
is zero for $e\gg 0$ if and only if $\dim R\ge 2$ and $\Spec R\setminus\{\m\}$
is connected in the Zariski topology.
\end{theorem}

We also obtain, by similar methods, the following theorem:

\begin{theorem}\label{components}
Let $(R,\m)$ be a complete local ring of positive dimension, with an
algebraically closed coefficient field of positive characteristic. Then the
number of connected components of $\Spec R\setminus\{\m\}$ is
$$
1+\dim_K \bigcap_{e\in\NN}F^e(H^1_\m(R)).
$$
\end{theorem}

In \S \ref{algorithm} we describe how this provides an algorithm to determine
the number of geometrically connected components of algebraic sets defined over
a finite field: computer algebra algorithms for primary decomposition can be
used to determine the number of connected components over finite extensions of
the fields $\FF_p$ or $\QQ$, but not over the algebraic closures of these
fields. In the case of characteristic zero, de~Rham cohomology allows for the
computation of the number of geometrically connected components via $D$-module
methods, \cite{W-proj}, and we show that the Frobenius provides analogous
methods in the case of positive characteristic.

Theorem~\ref{connected} is obtained in \cite{Ly-dual} as a corollary of the
following two theorems of Lyubeznik and Peskine-Szpiro:

\begin{theorem}\cite[Theorem~1.1]{Ly-dual}\label{dual}
Let $(A,\M)$ be a regular local ring containing a field of positive
characteristic, and let $\A$ be an ideal of $A$. Then $H^i_\A(A)=0$ if and only
if there exists an integer $e\ge 1$ such that the $e$-th Frobenius iteration
$$
F^e:H^{\dim A-i}_\M(A/\A) \to H^{\dim A-i}_\M(A/\A)
$$
is the zero map.
\end{theorem}

\begin{theorem}\cite[Chapter~III, Theorem~5.5]{PS}\label{thm-PS}
Let $(A,\M)$ be a complete regular local ring with a separably closed
coefficient field of positive characteristic, and let $\A$ be an ideal of $A$.
Then $H^i_\A(A)=0$ for $i\ge\dim A-1$ if and only if $\dim (A/\A)\ge 2$ and
$\Spec(A/\A)\setminus\{\M\}$ is connected.
\end{theorem}

Our proof of Theorem~\ref{connected} is \lq\lq simple\rq\rq\ in the sense that
it does not rely on vanishing theorems such as those of \cite{PS}---indeed, the
only ingredient, aside from elementary considerations, is the local duality
theorem.

Results analogous to Theorem~\ref{thm-PS} were proved by Hartshorne in the
projective case \cite[Theorem~7.5]{CDAV}, and by Ogus in equicharacteristic
zero using de~Rham cohomology \cite[Corollary~2.11]{Ogus}. Combining these
results, one has:

\begin{theorem}
Let $(A,\M)$ be a regular local ring containing a field, and let $\A$ be an
ideal of $A$. Then $H^i_\A(A)=0$ for $i\geq \dim A-1$ if and only if
\begin{enumerate}
\item $\dim(A/\A)\geq 2$, and
\item $\Spec(A/\A)\setminus\{\M\}$ is formally geometrically connected (see
Definition \ref{defn-stuff}).
\end{enumerate}
\end{theorem}

Huneke and Lyubeznik \cite[Theorem~2.9]{HL} gave a characteristic free proof of
this using a generalization of a result of Faltings, \cite[Satz~1]{Faltings}.
Some other applications of local cohomology theory which yield strong results
on the connectedness properties of algebraic varieties may be found in the
papers \cite{BR} and \cite{HH}, where the authors obtain generalizations of
Faltings' connectedness theorem.

\section{Preliminary remarks}

\noindent {\bf Notation:} When $R$ is the homomorphic image of a ring $A$, we
use upper-case letters $\P,\Q,\M,\A,\B$ for ideals of $A$, and corresponding
lower-case letters $\p,\q,\m,\a,\b$ for their images in $R$.

\begin{defn}\label{defn-stuff}
Let $(R,\m)$ be a local ring. A field $K\subseteq R$ is a \emph{coefficient
field} for $R$ if the composition $K\into R\twoheadrightarrow R/\m$ is an
isomorphism. Every complete local ring containing a field has a coefficient
field.

We recall some notions from \cite[Chapitre~VIII]{Raynaud}. Let $(R,\m,K)$ be a
local ring and let $\overline{f(T)}\in K[T]$ denote the image of a polynomial
$f(T)\in R[T]$. Then $R$ is \emph{Henselian} if for every monic polynomial
$f(T)\in R[T]$, every factorization of $\overline{f(T)}$ as a product of
relatively prime monic polynomials in $K[T]$ lifts to a factorization of $f(T)$
as a product of monic polynomials in $R[T]$. Hensel's Lemma is precisely the
statement that every complete local ring is Henselian. The \emph{Henselization}
of a local ring $R$ is a local ring $R^\mathrm{h}$, with the property that
every local homomorphism from $R$ to a Henselian local ring factors uniquely
through $R^\mathrm{h}$. The ring $R^\mathrm{h}$ is obtained by taking the
direct limit of all local \'etale extensions $S$ of $R$ for which $(R,\m)\to
(S,\n)$ induces an isomorphism of residue fields $R/\m\overset{\cong}\to S/\n$.

A local ring $(R,\m,K)$ is said to be \emph{strictly Henselian} if it is
Henselian and its residue field $K$ is separably closed. It is easily seen that
$R$ is strictly Henselian if and only if every monic polynomial $f(T)\in R[T]$
for which $\overline{f(T)} \in K[T]$ is separable splits into linear factors in
$R[T]$. Every local ring has a \emph{strict Henselization} $R^\mathrm{sh}$,
such that every local homomorphism from $R$ to a strictly Henselian ring
factors through $R^\mathrm{sh}$. The strict Henselization of a field $K$ is its
separable closure $K^\sep$. In general, the strict Henselization of a local
ring $(R,\m,K)$ is obtained by fixing an embedding $\iota:K\to K^\sep$, and
taking the direct limit of local \'etale extensions $(S,\n,L)$ of $(R,\m,K)$
with $L\into K^\sep$, for which the induced map $K\to L\to K^\sep$ agrees with
$\iota:K\to K^\sep$.

The \emph{punctured spectrum} of a local ring $(R,\m)$ is the set $\Spec
R\setminus\{\m\}$, with the topology induced by the Zariski topology on $\Spec
R$. We say that \emph{the punctured spectrum of $R$ is formally geometrically
connected} if the punctured spectrum of $\hat{\hat{R}^{\mathrm{sh}}}$, the
completion of the strict Henselization of the completion of $R$, is connected.
If $R$ is an $\NN$-graded ring which is finitely generated over a field
$R_0=K$, then $\Proj R$ is said to be \emph{geometrically connected} if
$\Proj(R\otimes_K K^\sep)$, is connected.
\end{defn}

\begin{defn}
Let $\a$ be an ideal of a ring $R$. A ring homomorphism $\varphi:R\to S$
induces a map of local cohomology modules $H^i_\a(R)\overset{\varphi}\to
H^i_{\a S}(S)$. In particular, if $R$ contains a field of characteristic $p>0$,
then the Frobenius homomorphism $F:R\to R$ induces an additive map
$$
H^i_\a(R)\overset{F}\to H^i_{\a^{[p]}}(R)=H^i_\a(R),
$$
called the \emph{Frobenius action} on $H^i_\a(R)$. An element $\eta\in
H^i_\a(R)$ is \emph{$F$-torsion} if there exists $e\in\NN$ such that
$F^e(\eta)=0$. The module $H^i_\a(R)$ is \emph{$F$-torsion} if every element of
$H^i_\a(R)$ is $F$-torsion. The image of $F^e$ is not, in general, an
$R$-module. However it is a $K$-vector space when $K$ is perfect, and in this
case the \emph{$F$-stable} part of $H^i_\a(R)$ is the $K$-vector space
$$
H^i_\a(R)_\st=\bigcap_{e\in\NN}F^e(H^i_\a(R)).
$$
Some results about $F$-torsion modules and $F$-stable subspaces are summarized
in \S\ref{appendix}. For a very general theory of $F$-modules, we refer the
reader to \cite{Ly-F}.
\end{defn}

\begin{remark}\label{graph}
Consider a local ring $(R,\m)$ of positive dimension. The punctured spectrum of
$R$ is disconnected if and only if the minimal primes of $R$ can be partitioned
into two sets $\p_1,\dots,\p_m$ and $\q_1,\dots,\q_n$ such that
$\rad(\p_i+\q_j)=\m$ for all pairs $\p_i,\q_j$. Consider the graph $\Gamma$
whose vertices are the minimal primes of $R$, and there is an edge between
minimal primes $\p$ and $\p'$ if and only if $\rad(\p+\p')\neq\m$. It follows
that the punctured spectrum of $R$ is connected if and only if the graph
$\Gamma$ is connected. If the graph $\Gamma$ is connected, take a
\emph{spanning tree} for it, i.e., a connected acyclic subgraph, containing all
the vertices of $\Gamma$. This spanning tree must contain a vertex $\p_i$ with
only one edge, so $\Gamma\setminus\{\p_i\}$ is connected as well.

Let $\P_1,\dots,\P_n$ be incomparable prime ideals of a local domain $A$. Then
their images $\p_1,\dots,\p_n$ are precisely the minimal primes of the ring
$R=A/(\P_1\cap\dots\cap\P_n)$. From the above discussion, we conclude that if
the punctured spectrum of $R$ is connected, then there exists $i$ such that the
punctured spectrum of the ring
$$
A/(\P_1\cap\dots\cap\hat{\P}_i\cap\dots\cap\P_n)
$$
is connected as well.
\end{remark}

Theorems~\ref{connected} and \ref{components} assert that connectedness issues
for $\Spec R\setminus\{\m\}$ are determined by the Frobenius action on
$H^1_\m(R)$. We next record an observation about the length of $H^1_\m(R)$.

\begin{prop}\label{length}
Let $(R,\m)$ be a local ring which is a homomorphic image of a Gorenstein
domain. Then $H^1_\m(R)$ has finite length if and only if $\ann_R\p=0$ for
every prime ideal $\p$ of $R$ with $\dim R/\p=1$.
\end{prop}

\begin{proof}
If $\dim R=0$, then $H^1_\m(R)=0$, and $R$ has no primes with $\dim R/\p=1$. If
$\dim R=1$, then $H^1_\m(R)$ has infinite length and $\dim R/\p=1$ for some
minimal, hence associated, prime $\p$ of $R$. For the rest of the proof we
hence assume that $\dim R\ge 2$.

Let $R=A/\Q$ where $A$ is a Gorenstein domain. Localizing $A$ at the inverse
image of $\m$, we may assume that $(A,\M)$ is a local ring. Using local duality
over $A$, the module $H^1_\m(R)=H^1_\M(A/\Q)$ has finite length if and only if
$\Ext^{\dim A-1}_A(A/\Q,A)$ has finite length as an $A$-module. Since
$\Ext^{\dim A-1}_A(A/\Q,A)$ is finitely generated, this is equivalent to the
vanishing of
$$
\Ext^{\dim A-1}_A(A/\Q,A)_\P=\Ext^{\dim A-1}_{A_\P}(A_\P/\Q A_\P,A_\P)
$$
for all $\P\in\Spec A\setminus\{\M\}$. Using local duality over the Gorenstein
local ring $(A_\P,\P A_\P)$, this is equivalent to the vanishing of
$$
H^{\dim A_\P-\dim A+1}_{\P A_\P}(A_\P/\Q A_\P) =H^{\dim A_\P-\dim A+1}_{\p
R_\p}(R_\p)
$$
for all $\P\in\Spec A\setminus\{\M\}$. This local cohomology module vanishes
for $\P\notin V(\Q)$. Since $\dim A_\P-\dim A+1\le 0$ for $\P\in\Spec
A\setminus\{\M\}$, we need only consider primes $\P\in V(\Q)$ with $\dim
A_\P=\dim A-1$. Since $A$ is a catenary local domain, $\dim A_\P$ equals $\dim
A-1$ precisely when $\dim A/\P=1$, which is equivalent to $\dim R/\p=1$. Hence
$H^1_\m(R)$ has finite length if and only if $H^0_{\p R_\p}(R_\p)=H^0_\p(R)$
vanishes for all $\p\in\Spec R$ with $\dim R/\p=1$, i.e., if and only if
$\ann_R\p=0$ for all $\p\in\Spec R$ with $\dim R/\p=1$.
\end{proof}

\section{Main results}\label{main}

\begin{theorem}\label{domain1}
Let $(R,\m)$ be a strictly Henselian local domain containing a field of
positive characteristic. If $R$ is a homomorphic image of a Gorenstein domain
and $\dim R\ge 2$, then $H^1_\m(R)$ is $F$-torsion.
\end{theorem}

\begin{proof}
Suppose there exists $\eta\in H^1_\m(R)$ which is not $F$-torsion. Since $R$ is
a domain, Proposition~\ref{length} implies that $H^1_\m(R)$ has finite length.
Hence for all integers $e\gg 0$, the element $F^e(\eta)$ belongs to the
$R$-module spanned by $\eta,F(\eta),F^2(\eta),\dots,F^{e-1}(\eta)$. Amongst all
equations of the form
\begin{eqnarray}
\label{eqn-min-eta} F^{e+k}(\eta)+r_1 F^{e+k-1}(\eta)+\dots+r_eF^k(\eta)=0
\end{eqnarray}
with $r_i\in R$ for all $i$, choose one where the number of nonzero
coefficients $r_i$ that occur is minimal. We claim that $r_e$ must be a unit.
Note that $H^1_\m(R)$ is killed by $\m^{q'}$ for some $q'=p^{e'}$. If
$r_e\in\m$, then applying $F^{e'}$ to equation (\ref{eqn-min-eta}), we get
$$
F^{e'+e+k}(\eta)+r_1^{q'} F^{e'+e+k-1}(\eta)+\dots+r_e^{q'} F^{e'+k}(\eta)=0.
$$
But $r_e^{q'}F^{e'+k}(\eta)\in\m^{q'}H^1_\m(R)=0$, so this is an equation with
fewer nonzero coefficients, contradicting the minimality assumption. This shows
that $r_e\in R$ is a unit. Since $\eta$ is not $F$-torsion, neither is
$F^k(\eta)$, so after replacing $\eta$ if necessary, we have an equation
\begin{eqnarray}
\label{eqn-other-eta} F^e(\eta)+r_1 F^{e-1}(\eta)+\dots+r_e\eta=0
\end{eqnarray}
where $r_e$ is a unit and $\eta\in H^1_\m(R)$ is not $F$-torsion. Let
$\eta=[(y_1/x_1,\dots,y_d/x_d)]$ where $H^1_\m(R)$ is regarded as the
cohomology of a \v Cech complex on a system of parameters $x_1,\dots,x_d$ for
$R$. Then (\ref{eqn-other-eta}) implies that there exists $r_{e+1}\in R$ such
that each $y_i/x_i\in R_{x_i}$ is a root of the polynomial
$$
f(T)=T^{p^e}+r_1 T^{p^{e-1}}+\dots+r_e T+r_{e+1}\in R[T].
$$
Now $f'(T)=r_e$ is a unit, so $\overline{f(T)}\in R/\m[T]$ is a separable
polynomial. Since $R$ is strictly Henselian, the polynomial $f(T)$ splits in
$R[T]$, and hence any root of $f(T)$ in the fraction field of $R$ must be an
element of $R$. In particular, $y_1/x_1=\dots=y_d/x_d\in R$, and so $\eta=0$.
\end{proof}

We next prove the connectedness criterion, Theorem~\ref{connected}. By
Proposition~\ref{torsionindex}, the module $H^1_\m(R)$ is $F$-torsion if and
only if there exists $e$ such that $F^e(H^1_\m(R))=0$. In view of this, the
following theorem is equivalent to Theorem~\ref{connected}.

\begin{theorem}\label{connected2}
Let $(R,\m)$ be a local ring with $\dim R>0$, which contains a field of
positive characteristic. Then $H^1_\m(R)$ is $F$-torsion if and only if $\dim
R\ge 2$ and the punctured spectrum of $R$ is formally geometrically connected.
\end{theorem}

\begin{proof}
Quite generally, for a local ring $(R,\m)$ we have $H^i_\m(\hat{R})=H^i_\m(R)$.
Moreover, $S=\hat{\hat{R}^\mathrm{sh}}$ is a faithfully flat extension of $R$,
and $H^i_\m(R)\otimes_R S\cong H^i_{\m S}(S)$ is $F$-torsion if and only if
$H^i_\m(R)$ is $F$-torsion. Hence we may assume that $R$ is a complete local
ring with a separably closed coefficient field.

Suppose that $H^1_\m(R)$ is $F$-torsion. The local cohomology module $H^{\dim
R}_\m(R)$ is not $F$-torsion by Lemma~\ref{list}, so $\dim R\ge 2$. Let $\a$
and $\b$ be ideals of $R$ such that $\a+\b$ is $\m$-primary and $\a\cap\b=0$.
Let
$$
x_1=y_1+z_1,\quad\dots,\quad x_d=y_d+z_d
$$
be a system of parameters for $R$ where $y_i\in\a$ and $z_i\in\b$. Since
$\a\b\subseteq\a\cap\b=0$, we have $y_iz_j=0$ for all $i,j$, and hence
$$
y_i(y_j+z_j)=y_j(y_i+z_i).
$$
These relations give an element of $H^1_\m(R)$ regarded as the cohomology of a
\v Cech complex on $x_1,\dots,x_d$, namely
$$
\eta=\left[\left(\frac{y_1}{x_1},\dots,\frac{y_d}{x_d}\right)\right] \in
H^1_\m(R).
$$
The hypothesis implies that $F^e(\eta)=0$ for some $e$, so there exists $q=p^e$
and $r\in R$ such that $(y_i/x_i)^q=r$ in $R_{x_i}$ for all $1\le i\le d$.
Hence there exists $t\in\NN$ such that $x_i^ty_i^q=rx_i^{q+t}$, i.e.,
$$
(y_i+z_i)^t y_i^q=r (y_i+z_i)^{q+t}.
$$
But $y_iz_i=0$, so these equations simplify to give
$(1-r)y_i^{q+t}=rz_i^{q+t}$. Since $R$ is a local ring, either $r$ or $1-r$
must be a unit. If $r$ is a unit, then $z_i^{q+t}\in\a$ for all $i$, and so
$\a$ is $\m$-primary. Similarly if $1-r$ is a unit, then $\b$ is $\m$-primary.
This proves that the punctured spectrum of $R$ is connected.

For the converse, assume that $\dim R\ge 2$ and that the punctured spectrum of
$R$ is connected. Let $\n$ denote the nilradical of $R$. Note that $\Spec R$ is
homeomorphic to $\Spec R/\n$. Moreover, $\n$ supports a Frobenius action and is
$F$-torsion. The long exact sequence of local cohomology relating $H^1_\m(R)$
and $H^1_\m(R/\n)$ implies that if $H^1_\m(R/\n)$ is $F$-torsion then so is
$H^1_\m(R)$, and hence there is no loss of generality in assuming that $R$ is
reduced. Let $R=A/(\P_1\cap\dots\cap\P_n)$ where $\P_1,\dots,\P_n$ are
incomparable prime ideals of a power series ring $A=K[[x_1,\dots,x_m]]$ over a
separably closed field $K$. We use induction on $n$ to prove that $H^1_\m(R)$
is $F$-torsion; the case $n=1$ follows from Theorem~\ref{domain1}, so we assume
$n>1$ below.

If $\dim R/\p_i=1$ for some $i$, then $\Spec R\setminus\{\m\}$ is the disjoint
union of $V(\p_i)\setminus\{\m\}$ and
$V(\p_1\cap\dots\cap\hat{\p}_i\cap\dots\cap\p_n)\setminus\{\m\}$, contradicting
the connectedness assumption. Hence $\dim R/\p_i\ge 2$ for all $i$. By
Remark~\ref{graph}, after relabeling the minimal primes if necessary, we may
assume that the punctured spectrum of $A/\Q$ is connected where
$\Q=\P_2\cap\dots\cap\P_n$. The short exact sequence
$$
0\to A/(\P_1\cap\Q)\to A/\P_1\oplus A/\Q\to A/(\P_1+\Q)\to 0
$$
induces a long exact sequence of local cohomology modules containing the piece
\begin{eqnarray}
\label{eqn-les} \quad H^0_\M(A/(\P_1+\Q))\to H^1_\M(A/(\P_1\cap\Q)) \to
H^1_\M(A/\P_1)\oplus H^1_\M(A/\Q).
\end{eqnarray}
Since $\rad(\P_1 +\P_i)\neq\M$ for some $i>1$, it follows that $\dim
A/(\P_1+\Q)\ge 1$. Proposition~\ref{list}~(1) now implies that
$H^0_\M(A/(\P_1+\Q))$ is $F$-torsion. By the inductive hypothesis,
$H^1_\M(A/\P_1)$ and $H^1_\M(A/\Q)$ are $F$-torsion as well. The exact sequence
(\ref{eqn-les}) implies that $H^1_\M(A/(\P_1\cap\Q))=H^1_\m(R)$ is $F$-torsion.
\end{proof}

The following lemma will be used in the proof of Theorem~\ref{components}.

\begin{lemma}\label{domain2}
Let $(R,\m)$ be a complete local domain with an algebraically closed
coefficient field of positive characteristic. Then $H^1_\m(R)_\st$, the
$F$-stable part of $H^1_\m(R)$, is zero.
\end{lemma}

\begin{proof}
If $\dim R=0$, then $H^1_\m(R)=0$, and if $\dim R\ge 2$, then the assertion
follows from Theorem~\ref{domain1}. The remaining case is $\dim R=1$.
Theorem~\ref{finiteness} implies that $H^1_\m(R)_\st$ has a vector space basis
$\eta_1,\dots,\eta_r$ such that $F(\eta_i)=\eta_i$.

Let $\eta\in H^1_\m(R)_\st$ be an element with $F(\eta)=\eta$. Considering
$H^1_\m(R)$ as the cohomology of a suitable \v Cech complex, let $\eta$ be the
class of $y/x$ in $R_x/R=H^1_\m(R)$, where $y\in R$ and $x\in\m$. Since
$F(\eta)=\eta$, there exists $r\in R$ such that
$$
\left(\frac{y}{x}\right)^p-\frac{y}{x}-r=0,
$$
and so $y/x\in R_x$ is a root of the polynomial $f(T)=T^p-T-r\in R[T]$. The
polynomial $\overline{f(T)}\in K[T]$ is separable and $R$ is strictly
Henselian, so $f(T)$ splits in $R[T]$. Since $y/x$ is a root of $f(T)$ in the
fraction field of $R$, it must then be an element of $R$, and hence $\eta=0$.
\end{proof}

\begin{proof}[Proof of Theorem~\ref{components}]
We may assume $R$ to be reduced by Lemma~\ref{nilradical}. First consider the
case where the punctured spectrum of $R$ is connected. If $\dim R\ge 2$, then
$H^1_\m(R)$ is $F$-torsion by Theorem~\ref{connected2}, so $H^1_\m(R)_\st=0$.
If $\dim R=1$, then $R$ is a domain, and Lemma~\ref{domain2} implies that
$H^1_\m(R)_\st=0$.

We continue by induction on the number of connected components of the punctured
spectrum of $R$. If the punctured spectrum of $R$ is disconnected, then
$R=A/(\A\cap\B)$ where $(A,\M)$ is a power series ring over the field $K$, and
$\A$ and $\B$ are radical ideals of $A$ which are not $\M$-primary, but $\A+\B$
is $\M$-primary. There is a short exact sequence
$$
0\to A/(\A\cap\B)\to A/\A\oplus A/\B\to A/(\A+\B)\to 0.
$$
Since $H^0_\M(A/\A)=H^0_\M(A/\B)=H^1_\M(A/(\A+\B))=0$, the resulting long exact
sequence of local cohomology modules gives a short exact sequence
$$
0\to H^0_\M(A/(\A+\B))\to H^1_\M(A/(\A\cap\B)) \to H^1_\M(A/\A)\oplus
H^1_\M(A/\B)\to 0.
$$
By Theorem~\ref{exact}, we have a $K$-vector space isomorphism
$$
H^1_\m(R)_\st=H^1_\M(A/(\A\cap\B))_\st\cong H^0_\M(A/(\A+\B))_\st \oplus
H^1_\M(A/\A)_\st\oplus H^1_\M(A/\B)_\st.
$$
Since $H^0_\M(A/(\A+\B))_\st =K$ by Proposition~\ref{list}~(3), the inductive
hypothesis completes the proof.
\end{proof}

We next record the graded versions of the results proved in this section:

\begin{theorem}\label{graded}
Let $R$ be an $\NN$-graded ring of positive dimension, which is finitely
generated over a field $R_0=K$ of characteristic $p>0$.
\begin{enumerate}
\item If $R$ is a domain with $\dim R\ge 2$, and $K$ is separably closed, then
$H^1_\m(R)$ is $F$-torsion.
\item The module $H^1_\m(R)$ is $F$-torsion if and only if $\dim R\ge 2$ and
$\Proj R$ is geometrically connected.
\item Let $K$ be a perfect field, and let $\overline{K}$ denote its algebraic
closure. Then the number of connected components of $\Proj(R\otimes_K
\overline{K})$ is
$$
1+\dim_K H^1_\m(R)_\st=1+\dim_K ([H^1_\m(R)]_0)_\st.
$$
\end{enumerate}
\end{theorem}

\begin{proof}
(1) Note that $H^1_\m(R)$ is a $\ZZ$-graded $R$-module, and that
$$
\qquad F:[H^1_\m(R)]_n\to [H^1_\m(R)]_{np}\quad \text{for all }n\in\ZZ.
$$
The module $H^1_\m(R)$ has finite length, so all elements of $H^1_\m(R)$ of
positive or negative degree are $F$-torsion; it remains to show that elements
$\eta\in[H^1_\m(R)]_0$ are $F$-torsion as well. Let $\eta$ be a element of
$[H^1_\m(R)]_0$ which is not $F$-torsion. As in the proof of
Theorem~\ref{domain1}, after a change of notation we may assume that
$$
F^e(\eta)+r_1 F^{e-1}(\eta)+\dots+r_e\eta=0
$$
where all $r_i$ are in $[R]_0=K$, and $r_e$ is nonzero. Let
$\eta=[(y_1/x_1,\dots,y_d/x_d)]$ where $H^1_\m(R)$ is regarded as the
cohomology of a homogeneous \v Cech complex. Then there exists $r_{e+1}\in K$
such that $y_i/x_i\in R_{x_i}$ is a root of the polynomial
$$
f(T)=T^{p^e}+r_1 T^{p^{e-1}}+\dots+r_e T+r_{e+1}\in K[T].
$$
But $f(T)$ is a separable polynomial, so it splits in $K[T]$. The element
$y_i/x_i=y_j/x_j$ is a root of $f(T)$ in the fraction field of $R$, so it must
be one of the roots of $f(T)$ in $K$. It follows that $\eta=0$, which completes
the proof of (1).

The proof of (2) is now similar to that of Theorem~\ref{connected2}, and is
left to the reader. For (3), note that $F^e(H^1_\m(R))$ is a $K$-vector space
since $K$ is perfect, and that
$$
\dim_K H^1_\m(R)_\st=\dim_{\overline{K}}H^1_\m(R\otimes_K \overline{K})_\st.
$$
Thus we may assume $K$ is algebraically closed, and the proof is then similar
to that of Theorem~\ref{components}.
\end{proof}

\begin{remark}
Theorem~\ref{graded}~(3) generalizes, in the case of positive characteristic,
the well-known fact that the number of connected components of $X=\Proj R$ is
$$
\dim_K H^0(X,\calO_X)=1+\dim_K [H^1_\m(R)]_0,
$$
where $R$ is an $\NN$-graded \emph{reduced} ring of positive dimension, which
is finitely generated over an algebraically closed field $R_0=K$. The point is
that in this case the Frobenius is bijective on $[H^1_\m(R)]_0$. To see this,
let
$$
\eta=\left[\left(\frac{y_1}{x_1},\dots,\frac{y_d}{x_d}\right)\right]\in
[H^1_\m(R)]_0
$$
be an element with $F(\eta)=0$, where $H^1_\m(R)$ is computed as the cohomology
of a suitable \v Cech complex. Then there exists a homogeneous element $r\in R$
with $(y_i/x_i)^p=r$ in $R_{x_i}$ for all $1\le i\le d$. Such an element $r$
must have degree zero, and hence must be an element of $K$. But then
$r^{1/p}\in K$, and, since $R$ is reduced, $y_i/x_i=r^{1/p}$ for all $i$. It
follows that
$$
\eta=[(r^{1/p},\dots,r^{1/p})]=0.
$$
To complete the argument, note that $[H^1_\m(R)]_0$ is a finite dimensional
$K$-vector space, and that if $\eta_1,\dots,\eta_n\in [H^1_\m(R)]_0$ are
linearly independent, then so are $F(\eta_1),\dots,F(\eta_n)$. It follows that
$F:[H^1_\m(R)]_0\to [H^1_\m(R)]_0$ is surjective.
\end{remark}

\section{$F$-purity}

A ring homomorphism $\varphi:R\to S$ is \emph{pure} if $\varphi\otimes 1
:R\otimes_R M\to S\otimes_R M$ is injective for every $R$-module $M$. If $R$ is
a ring containing a field of characteristic $p>0$, then $R$ is \emph{ $F$-pure}
if the Frobenius homomorphism $F:R\to R$ is pure. The notion was introduced by
Hochster and Roberts in the course of their study of rings of invariants in
\cite{HR1,HR2}.

Examples of $F$-pure rings include regular rings of positive characteristic and
their pure subrings. If $\a$ is generated by square-free monomials in the
variables $x_1,\dots,x_n$ and $K$ is a field of positive characteristic, then
$K[x_1,\dots,x_n]/\a$ is $F$-pure.

Goto and Watanabe classified one-dimensional $F$-pure rings in \cite{GW}: let
$(R,\m)$ be a local ring containing a field of positive characteristic such
that $R/\m=K$ is algebraically closed, $F:R\to R$ is finite, and $\dim R=1$.
Then $R$ is $F$-pure if and only if
$$
\hat{R}\cong K[[x_1,\dots,x_n]]/(x_ix_j:i<j).
$$

Two-dimensional $F$-pure rings have attracted a lot of attention: in
\cite{Wat-Fpure} Watanabe proved that $F$-pure normal Gorenstein local rings of
dimension two are either rational double points, simple elliptic singularities,
or cusp singularities. Watanabe also obtained a classification of
two-dimensional normal $\NN$-graded rings $R$ over an algebraically closed
field $R_0$, in terms of the associated $\QQ$-divisor on the curve $\Proj R$,
\cite{Wat-dim2}. In \cite{MS} Mehta and Srinivas obtained a classification of
two-dimensional $F$-pure normal singularities in terms of the resolution of the
singularity. Hara completed the classification of two-dimensional normal
$F$-pure singularities in terms of the dual graph of the minimal resolution of
the singularity, \cite{Hara}.

The results of \S\ref{main} imply that over separably closed fields, $F$-pure
domains of dimension two are Cohen-Macaulay. The point is that if $R$ is an
$F$-pure ring, then the Frobenius action $F:H^i_\m(R)\to H^i_\m(R)$ is an
injective map.

\begin{cor}\label{Fpure}
Let $R$ be a local ring with $\dim R\ge 2$, which contains a field of positive
characteristic. If $R$ is $F$-pure and the punctured spectrum of $R$ is
formally geometrically connected, then $\depth R\ge 2$.

In particular, if $R$ is a complete local $F$-pure domain of dimension two,
with a separably closed coefficient field, then $R$ is Cohen-Macaulay.
\end{cor}

\begin{proof}
An $F$-pure ring is reduced, so $H^0_\m(R)=0$. By Theorem~\ref{domain1},
$H^1_\m(R)$ is $F$-torsion. Since $R$ is $F$-pure, it follows that
$H^1_\m(R)=0$.
\end{proof}

In the graded case, we similarly have:

\begin{cor}
Let $R$ be an $\NN$-graded ring with $\dim R\ge 2$, which is finitely generated
over a field $R_0=K$ of positive characteristic. If $R$ is $F$-pure and $\Proj
R$ is geometrically connected, then $\depth R\ge 2$.
\end{cor}

In the following example, $R$ is a graded $F$-pure domain of dimension $2$, but
$\depth R=1$. The issue is that $\Proj R$ is connected though not geometrically
connected.

\begin{example}
Let $K$ be a field of characteristic $p> 2$, and $a\in K$ an element such that
$\sqrt{a}\notin K$. Let $R=K[x,y,x\sqrt{a},y\sqrt{a}]$. The domain $R$ has a
presentation
$$
R=K[x,y,u,v]/(u^2-ax^2, v^2-ay^2, uv-axy, vx-uy),
$$
and if $K^\sep$ denotes the separable closure of $K$, then
$$
R\otimes_K K^\sep\cong
K^\sep[x,y,u,v]/(u-x\sqrt{a},v-y\sqrt{a})(u+x\sqrt{a},v+y\sqrt{a}).
$$
Using a change of variables, $R\otimes_K K^\sep\cong
K^\sep[x',y',u',v']/(x',y')(u',v')$. Since $(x',y')(u',v')$ is a square-free
monomial ideal, $R\otimes_K K^\sep$ is $F$-pure and it follows that $R$ is
$F$-pure. However, $R$ is not Cohen-Macaulay since $x,y$ is a homogeneous
system of parameters with a non-trivial relation
$$
(x\sqrt{a})y=(y\sqrt{a})x.
$$
Using the \v Cech complex on $x,y$ to compute $H^1_\m(R)$, we see that it is a
$1$-dimensional $K$-vector space generated by the element
$$
\eta=\left[\left(\frac{x\sqrt{a}}{x},\frac{y\sqrt{a}}{y}\right)\right] \in
H^1_\m(R)
$$
corresponding to the relation above. Given $e\in\NN$, let $p^e=2k+1$. Then
$$
F^e(\eta)=a^k\eta,
$$
which is a nonzero element of $H^1_\m(R)$. Consequently $H^1_\m(R)$ is not
$F$-torsion, corresponding to that fact that $\Proj R$ is not geometrically
connected.
\end{example}

The corollaries obtained in this section imply that over a separably closed
field, a graded or complete local $F$-pure domain of dimension two is
Cohen-Macaulay. We record an example which shows that this is not true for
rings of higher dimension.

\begin{example}
Let $K$ be a field of characteristic $p>0$, and take
$$
A=K[x_1,\dots,x_d]/(x_1^d+\dots+x_d^d)
$$
where $d\ge 3$. Let $R$ be the Segre product of $A$ and the polynomial ring
$B=K[s,t]$. Then $\dim R=d$, and the K\"unneth formula for local cohomology
implies that
$$
H^{d-1}_{\m_R}(R)\ \cong\ [H^{d-1}_{\m_A}(A)]_0 \otimes_K [B]_0\ \cong\ K,
$$
so $R$ is not Cohen-Macaulay. If $p\equiv 1\mod d$ then $A$ is $F$-pure by
\cite[Proposition~5.21]{HR2}, hence $A\otimes_KB$ and its direct summand $R$
are $F$-pure as well.
\end{example}

\section{Algorithmic aspects}\label{algorithm}

Let $R$ be an $\NN$-graded ring, which is finitely generated over a finite
field $R_0=K$. We wish to determine the number of geometrically connected
components of the scheme $\Proj R$, i.e., the number of connected components of
$\Proj(R\otimes_K \overline{K})$, or, equivalently, of $\Proj(R\otimes_K
K^\sep)$. While primary decomposition algorithms such as those of \cite{EHV},
\cite{GTZ}, or \cite{SY}, may be used to determine the connected components of
$\Proj R$, there is computationally no hope of ``determining'' the connected
components over the algebraic closure, $\overline{K}$. However, simply finding
their number is much easier: by Theorem~\ref{graded}, this is
$1+\dim_K([H^1_\m(R)]_0)_\st$. Computing this number involves three steps:
\begin{enumerate}
\item Finding a good presentation of $[H^1_\m(R)]_0$;
\item Determining the Frobenius action on $[H^1_\m(R)]_0$ in terms of this
presentation;
\item Computing the dimension of the $F$-stable part, $([H^1_\m(R)]_0)_\st$.
\end{enumerate}

If $R=A/\A$ for a polynomial ring $A$, we first replace $\A$ by an ideal that
has the same radical as $\A$, but does not have the homogeneous maximal ideal
$\M$ as associated prime. This can be done by saturating $\A$ with respect to
$\M$; if desired, one may simply compute the radical of $\A$, but this is often
computationally expensive. Now, since $\M$ is not associated to $\A$, one can
find a homogeneous system of parameters $x_1,\dots,x_d$ for $R$ such that each
$x_i$ is a nonzerodivisor on $R$.

The length $\ell$ of $[H^1_\m(R)]_0$ may be computed by computing the length of
its graded dual $[\Ext^{n-1}_A(R,A(-n))]_0$, where $\dim A=n$. Of course, if
this length is zero, then
$X_{\overline{K}}$ is connected. Consider the Koszul cohomology modules
$$
H^1(x_1^t,\dots,x_d^t;R)=\frac{\left\{(a_1,\dots,a_d)\in R^d:a_ix_j^t=a_jx_i^t
\text{ for all }i<j\right\}}{\left\{(rx_1^t,\dots,rx_d^t):r\in R\right\}}.
$$
These modules have an $\NN$-grading, where for homogeneous elements $a_i\in R$,
we define the degree of $[(a_1,\dots,a_d)]\in H^1(x_1^t,\dots,x_d^t;R)$ as
$$
\deg[(a_1,\dots,a_d)]=\deg a_i-\deg x_i^t,
$$
which is independent of $i$. This ensures that for each $t$, the map
\begin{align*}
H^1(x_1^t,\dots,x_d^t;R) &\to H^1(x_1^{t+1},\dots,x_d^{t+1};R)\\
[(a_1,\dots,a_d)] &\longmapsto [(a_1x_1,\dots,a_dx_t)]
\end{align*}
preserves degrees. The module $H^1_\m(R)$ is the direct limit of these Koszul
cohomology modules, and the assumption that the $x_i$ are nonzerodivisors
ensures that the maps in the direct limit system are injective. The modules
$H^1(x_1^t,\dots,x_d^t;R)$ may be computed for increasing values of $t$, until
we arrive at an integer $N$ such that
$$
\ell\left([H^1(x_1^N,\dots,x_d^N;R)]_0\right)=\ell.
$$
This gives us a presentation for $[H^1_\m(R)]_0=[H^1(x_1^N,\dots,x_d^N;R)]_0$,
in terms of which we now analyze the Frobenius map. Replacing the $x_i$ by
their powers if needed, assume that $N=1$. Let
$$
\alpha=[(a_1,\dots,a_d)]\in [H^1(x_1,\dots,x_d;R)]_0,
$$
in which case, $F(\alpha)=[(a_1^p,\dots,a_d^p)]\in
[H^1(x_1^p,\dots,x_d^p;R)]_0$. Since the map
$$
[H^1(x_1,\dots,x_d;R)]_0\to [H^1(x_1^p,\dots,x_d^p;R)]_0
$$
coming from the direct limit system is bijective, it follows that $a_i^p\in
x_i^{p-1}R$ for each $1\le i\le d$. Setting $b_i=a_i^p/x_i^{p-1}$, we arrive at
$$
F(\alpha)=[(b_1,\dots,b_d)]\in [H^1(x_1,\dots,x_d;R)]_0.
$$
Using this description of Frobenius action on the finite dimensional $K$-vector
space $[H^1_\m(R)]_0=[H^1(x_1,\dots,x_d;R)]_0$, it is now straightforward to
compute the ranks of the vector spaces
$$
[H^1_\m(R)]_0\supseteq F([H^1_\m(R)]_0)\supseteq
F^2([H^1_\m(R)]_0)\supseteq\dots,
$$
and hence of the $F$-stable part, $([H^1_\m(R)]_0)_\st$.

\section{Appendix: $F$-torsion modules and $F$-stable vector
spaces}\label{appendix}

Let $R$ be a commutative ring containing a field $K$ of characteristic $p>0$. A
\emph{Frobenius action} on an $R$-module $M$ is an additive map $F:M\to M$ such
that $F(rm)=r^pF(m)$ for all $r\in R$ and $m\in M$. In this case $\ker F$ is a
submodule of $M$, and we have an ascending sequence of submodules of $M$,
$$
\ker F\subseteq\ker F^2 \subseteq\ker F^3\subseteq\dots.
$$
The union of these is the \emph{$F$-nilpotent} submodule of $M$, denoted
$M_{\nil}=\bigcup_{e\in\NN}\ker F^e$. We say $M$ is \emph{$F$-torsion} if
$M_{\nil}=M$.

\begin{prop}\label{torsionindex}
Let $(R,\m)$ be a local ring containing a field of positive characteristic, and
let $M$ be an Artinian $R$-module with a Frobenius action. Then there exists
$e\in\NN$ such that $F^e(M_{\nil})=0$.

In particular, an Artinian module $M$ is $F$-torsion if and only if $F^e(M)=0$
for some $e\in\NN$.
\end{prop}

\begin{proof}
This is proved in \cite[Proposition~1.11]{HS} under the hypothesis that $R$ is
a complete local ring with a perfect coefficient field. The general case may be
concluded from this, but a more elegant approach is via Lyubeznik's theory of
$F$-modules, see \cite[Proposition~4.4]{Ly-F}.
\end{proof}

If $R$ is a ring containing a perfect field $K$ of positive characteristic and
$M$ is an $R$-module with a Frobenius action, then $F(M)$ is a $K$-vector
space, and we have a descending sequence of $K$-vector spaces
$$
F(M)\supseteq F^2(M)\supseteq F^3(M)\supseteq\dots.
$$
The \emph{$F$-stable} part of $M$ is the vector space
$M_\st=\bigcap_{e\in\NN}F^e(M)$.

\begin{prop}\label{list}
Let $(R,\m,K)$ be a local ring of dimension $d$ which contains a field of
positive characteristic.
\begin{enumerate}
\item $H^0_\m(R)$ is $F$-torsion if and only if $d>0$.
\item $H^d_\m(R)$ is not $F$-torsion. \item If $d=0$ and $K$ is perfect, then
$H^0_\m(R)_\st=R_\st=K$.
\end{enumerate}
\end{prop}

\begin{proof}
(1) If $d=0$ then $H^0_\m(R)=R$, which is not $F$-torsion. If $d>0$ then
$H^0_\m(R)$ is contained in $\m$. Since every element of $H^0_\m(R)$ is killed
by a power of $\m$, it follows that each element is nilpotent. (See also
\cite[Corollary~4.6(a)]{Ly-dual}.)

(2) View $H^d_\m(R)$ as the cohomology of a \v Cech complex on a system of
parameters $x_1,\dots,x_d$ for $R$, and let $\eta=[1+(x_1,\dots,x_d)]\in
H^d_\m(R)$. For all $e_0\in\NN$, the collection of elements $F^e(\eta)$ with
$e>e_0$ generates $H^d_\m(R)$ as an $R$-module. Hence $F^{e_0}(\eta)$ cannot be
zero by Grothendieck's nonvanishing theorem.

(3) Since $\m$ is nilpotent in this case, for integers $e\gg 0$ we have
$$
F^e(H^0_\m(R))=F^e(R)=\{x^{p^e}:x\in R\}=\{(y+z)^{p^e}:y\in K, z\in\m\}=K.
$$
\end{proof}

\begin{theorem}\label{finiteness}
Let $(R,\m)$ be a local ring with a perfect coefficient field $K$ of positive
characteristic. Let $M$ be an Artinian $R$-module with a Frobenius action. Then
$M_\st$ is a finite dimensional $K$-vector space, and $F:M_\st\to M_\st$ is an
automorphism of the Abelian group $M_\st$.

If $K$ is algebraically closed, then there exists a $K$-basis $e_1,\dots,e_n$
for $M_\st$ such that $F(e_i)=e_i$ for all $1\le i\le n$.
\end{theorem}

\begin{proof}
For the finiteness assertion, see \cite[Theorem~1.12]{HS} or
\cite[Proposition~4.9]{Ly-F}. It is easily seen that $F:M_\st\to M_\st$ is an
automorphism whenever $M_\st$ is finite dimensional. The existence of the
special basis when $K$ is algebraically closed follows from
\cite[Proposition~5,~page~233]{Dieudonne}.
\end{proof}

\begin{theorem}\cite[Theorem~1.13]{HS}\label{exact}
Let $(R,\m)$ be a complete local ring with an algebraically closed coefficient
field of positive characteristic. Let $L,M,N$ be $R$-modules with Frobenius
actions, such that we have a commutative diagram
$$
\begin{CD}
0@>>>L@>{\alpha}>>M@>{\beta}>>N@>>>0\\
@. @VFVV @VFVV @VFVV @.\\
0@>>>L@>{\alpha}>>M@>{\beta}>>N@>>>0
\end{CD}
$$
with exact rows. If $L$ is Noetherian and $N$ is Artinian, then the $F$-stable
parts form a short exact sequence
$$
0\to L_\st\to M_\st\to N_\st \to 0.
$$
\end{theorem}

\begin{prop}\label{nilradical}
Let $(R,\m,K)$ be a complete local ring with an algebraically closed
coefficient field of positive characteristic. Let $\n$ denote the nilradical of
$R$. Then for all $i\ge 0$, the natural map $H^i_\m(R)\to H^i_\m(R/\n)$, when
restricted to $F$-stable subspaces, gives an isomorphism
$$
H^i_\m(R)_\st\overset{\cong}\to H^i_\m(R/\n)_\st.
$$
\end{prop}

\begin{proof}
Let $k$ an integer such that $\n^{p^k}=0$. The short exact sequence
$$
0\to\n\to R\to R/\n\to 0
$$
induces a long exact sequence of local cohomology modules
$$
\cdots\to H^i_\m(\n)\overset{\alpha}\to H^i_\m(R) \overset{\beta}\to
H^i_\m(R/\n)\overset{\gamma}\to H^{i+1}_\m(\n)\to\cdots.
$$
Consider an element $\mu\in\ker(\beta)\cap H^i_\m(R)_\st$. Then
$\mu\in\image(\alpha)$, so $F^k(\mu)=0$. The Frobenius action on
$H^i_\m(R)_\st$ is an automorphism, so $\mu=0$, and hence the map
$H^i_\m(R)_\st\to H^i_\m(R/\n)_\st$ is injective.

To complete the proof it suffices, by Theorem~\ref{finiteness}, to consider an
element $\eta\in H^i_\m(R/\n)_\st$ with $F(\eta)=\eta$, and prove that it lies
in the image of $H^i_\m(R)_\st$. Now $\gamma(\eta)\in H^{i+1}_\m(\n)$ so
$F^k(\gamma(\eta))=0$, and therefore $F^k(\eta)=\eta\in\ker(\gamma)$.

Let $\eta=\beta(\mu)$ for some element $\mu\in H^i_\m(R)$. Then
$\beta(F(\mu)-\mu)=0$, which implies that $F(\mu)-\mu \in \image(\alpha)$.
Consequently $F^k(F(\mu)-\mu)=0$, which shows that $F^{k+1}(\mu)=F^k(\mu)$, and
hence that $F^k(\mu)\in H^i_\m(R)_\st$. Since
$$
\beta(F^k(\mu))=F^k(\beta(\mu))=F^k(\eta)=\eta,
$$
we are done.
\end{proof}


\begin{thebibliography}{Wat2}

\bibitem[BR]{BR}
M. Brodmann and J. Rung, \emph{Local cohomology and the connectedness dimension
in algebraic varieties}, Comment. Math. Helv. {\bf 61} (1986), 481--490.

\bibitem[Di]{Dieudonne}
J. Dieudonn\'e, \emph{Lie groups and Lie hyperalgebras over a field of
characteristic $p>0$. II}, Amer. J. Math. {\bf 77} (1955), 218--244.

\bibitem[EHV]{EHV}
D. Eisenbud, C. Huneke, and W. Vasconcelos, \emph{Direct methods for primary
decomposition}, Invent. Math. {\bf 110} (1992), 207--235.

\bibitem[Fa]{Faltings}
G. Faltings, \emph{\"Uber lokale Kohomologiegruppen hoher Ordnung}, J. Reine
Angew. Math. {\bf 313} (1980), 43--51.

\bibitem[GTZ]{GTZ}
P. Gianni, B. Trager, and G. Zacharias, \emph{Gr\"obner bases and primary
decomposition of polynomial ideals}, J. Symbolic Comput. {\bf 6} (1988),
149--167.

\bibitem[GW]{GW}
S. Goto and K.-i. Watanabe, \emph{The structure of one-dimensional $F$-pure
rings}, J. Algebra {\bf 49} (1977), 415--421.

\bibitem[HaN]{Hara}
N. Hara, \emph{Classification of two-dimensional $F$-regular and $F$-pure
singularities} Adv. Math. {\bf 133} (1998), 33--53.

\bibitem[HaR]{CDAV}
R. Hartshorne, \emph{Cohomological dimension of algebraic varieties}, Ann. of
Math. (2) {\bf 88} (1968), 403--450.

\bibitem[HS]{HS}
R. Hartshorne and R. Speiser, \emph{Local cohomological dimension in
characteristic $p$}, Ann. of Math. (2) {\bf 105} (1977), 45--79.

\bibitem[HH]{HH}
M. Hochster and C. Huneke, \emph{Indecomposable canonical modules and
connectedness}, in: Commutative algebra: syzygies, multiplicities, and
birational algebra (South Hadley, MA, 1992), Contemp. Math. {\bf 159}, Amer.
Math. Soc., Providence, RI, 1994, 197--208.

\bibitem[HR1]{HR1}
M. Hochster and J. Roberts, \emph{Rings of invariants of reductive groups
acting on regular rings are Cohen-Macaulay}, Advances in Math. {\bf 13} (1974),
115--175.

\bibitem[HR2]{HR2}
M. Hochster and J. Roberts, \emph{The purity of the Frobenius and local
cohomology}, Advances in Math. {\bf 21} (1976), 117--172.

\bibitem[HL]{HL}
C. Huneke and G. Lyubeznik, \emph{On the vanishing of local cohomology
modules}, Invent. Math. {\bf 102} (1990), 73--93.

\bibitem[Ly1]{Ly-F}
G. Lyubeznik, \emph{$F$-modules: applications to local cohomology and
$D$-modules in characteristic $p>0$}, J. Reine Angew. Math. {\bf 491} (1997),
65--130.

\bibitem[Ly2]{Ly-dual}
G. Lyubeznik, \emph{On the vanishing of local cohomology in characteristic
$p>0$}, Compositio Math. {\bf 142} (2006), 207--221.

\bibitem[MS]{MS}
V. B. Mehta and V. Srinivas, \emph{Normal $F$-pure surface singularities}, J.
Algebra {\bf 143} (1991), 130--143.

\bibitem[Og]{Ogus}
A. Ogus, \emph{Local cohomological dimension of algebraic varieties}, Ann. of
Math. (2) {\bf 98} (1973), 327--365.

\bibitem[PS]{PS}
C. Peskine and L. Szpiro, \emph{Dimension projective finie et cohomologie
locale}, Inst. Hautes \'Etudes Sci. Publ. Math. {\bf 42} (1973), 47--119.

\bibitem[Ra]{Raynaud}
M. Raynaud, \emph{Anneaux locaux hens\'eliens}, Lecture Notes in Mathematics
{\bf 169}, Springer-Verlag, Berlin-New York, 1970.

\bibitem[SY]{SY}
T. Shimoyama and K. Yokoyama, \emph{Localization and primary decomposition of
polynomial ideals}, J. Symbolic Comput. {\bf 22} (1996), 247--277.

\bibitem[Wal]{W-proj}
U. Walther, \emph{Algorithmic determination of the rational cohomology of
complex varieties via differential forms}, Contemp. Math. {\bf 286} (2001),
185--206.

\bibitem[Wat1]{Wat-Fpure}
K.-i. Watanabe, \emph{Study of $F$-purity in dimension two}, in: Algebraic
Geometry and Commutative Algebra in honor of Masayoshi Nagata, Vol. II,
Kinokuniya, Tokyo, 1988, 791--800.

\bibitem[Wat2]{Wat-dim2}
K.-i. Watanabe, \emph{$F$-regular and $F$-pure normal graded rings}, J. Pure
Appl. Algebra {\bf 71} (1991), 341--350.

\end{thebibliography}
\end{document}